\documentclass{elsart3-1}
\usepackage{amssymb}
\usepackage[english,francais]{babel}

\newtheorem{theorem}{Theorem}[section]
\newtheorem{lemma}[theorem]{Lemma}
\newtheorem{proposition}[theorem]{Proposition}
\newtheorem{corollary}[theorem]{Corollary}
\newtheorem{definition}[theorem]{Definition\rm}
\newtheorem{remark}{\it Remark\/}
\newtheorem{example}{\it Example\/}
\renewcommand{\theequation}{\arabic{equation}}
\def\og{\leavevmode\raise.3ex\hbox{$\scriptscriptstyle\langle\!\langle$~}}
\def\fg{\leavevmode\raise.3ex\hbox{~$\!\scriptscriptstyle\,\rangle\!\rangle$}}

\newcommand{\nc}[2]{\newcommand{#1}{#2}}
\newcommand{\rnc}[2]{\renewcommand{#1}{#2}}
\rnc{\theequation}{\thesection.\arabic{equation}}
\def\note#1{{}}
\def\Label#1{\label{#1}\ifmmode\llap{[#1] }\else
\marginpar{\smash{\hbox{[#1]}}}\fi}

\nc{\beq}{\begin{equation}}
\nc{\eeq}{\end{equation}}
\rnc{\[}{\beq}
\rnc{\]}{\eeq}
\rnc{\Label}{\label}
\nc{\ba}{\begin{array}}
\nc{\ea}{\end{array}}
\nc{\bea}{\begin{eqnarray}}
\nc{\beas}{\begin{eqnarray*}}
\nc{\eeas}{\end{eqnarray*}}
\nc{\eea}{\end{eqnarray}}
\nc{\be}{\begin{enumerate}}
\nc{\ee}{\end{enumerate}}
\nc{\bd}{\begin{diagram}}
\nc{\ed}{\end{diagram}}
\nc{\bi}{\begin{itemize}}
\nc{\ei}{\end{itemize}}
\nc{\bpr}{\begin{proposition}}

\nc{\bth}{\begin{theorem}}
\nc{\ble}{\begin{lemma}}
\nc{\bco}{\begin{corollary}}
\nc{\bre}{\begin{remark}}
\nc{\bex}{\begin{example}}
\nc{\bde}{\begin{definition}}
\nc{\ede}{\end{definition}}
\nc{\epr}{\end{proposition}}
\nc{\ethe}{\end{theorem}}
\nc{\ele}{\end{lemma}}
\nc{\eco}{\end{corollary}}
\nc{\ere}{\hfill\mbox{$\losenge$}\end{remark}}
\nc{\eex}{\hfill\mbox{$\losenge$}\end{example}}
\nc{\bpf}{{\it Proof.~~}}
\nc{\epf}{\hfill\mbox{$\square$}\vspace*{3mm}}
\nc{\hsp}{\hspace*}
\nc{\vsp}{\vspace*}

\def\ot{\otimes}
\nc{\te}{\!\ot\!}
\nc{\bmlp}{\mbox{\boldmath$\left(\right.$}}
\nc{\bmrp}{\mbox{\boldmath$\left.\right)$}}
\nc{\LAblp}{\mbox{\LARGE\boldmath$($}}
\nc{\LAbrp}{\mbox{\LARGE\boldmath$)$}}
\nc{\Lblp}{\mbox{\Large\boldmath$($}}
\nc{\Lbrp}{\mbox{\Large\boldmath$)$}}
\nc{\lblp}{\mbox{\large\boldmath$($}}
\nc{\lbrp}{\mbox{\large\boldmath$)$}}
\nc{\blp}{\mbox{\boldmath$($}}
\nc{\brp}{\mbox{\boldmath$)$}}
\nc{\LAlp}{\mbox{\LARGE $($}}
\nc{\LArp}{\mbox{\LARGE $)$}}
\nc{\Llp}{\mbox{\Large $($}}
\nc{\Lrp}{\mbox{\Large $)$}}
\nc{\llp}{\mbox{\large $($}}
\nc{\lrp}{\mbox{\large $)$}}
\nc{\lbc}{\mbox{\Large\boldmath$,$}}
\nc{\lc}{\mbox{\Large$,$}}
\nc{\Lall}{\mbox{\Large$\forall\;$}}
\nc{\bc}{\mbox{\boldmath$,$}}
\nc{\ra}{\rightarrow}
\nc{\ci}{\circ}
\nc{\cc}{\!\ci\!}
\nc{\lra}{\longrightarrow}
\nc{\imp}{\Rightarrow}
\rnc{\iff}{\Leftrightarrow}
\nc{\inc}{\mbox{$\,\subseteq\;$}}
\rnc{\subset}{\inc}
\nc{\0}{\sb{(0)}}
\nc{\1}{\sb{(1)}}
\nc{\2}{\sb{(2)}}
\nc{\3}{\sp{(3)}}
\nc{\4}{\sp{(4)}}
\nc{\5}{\sp{(5)}}
\nc{\6}{\sp{(6)}}
\nc{\7}{\sp{(7)}}

\def\<{\langle}
\def\>{\rangle}

\def\id{\mbox{$\mathop{\mbox{\rm id}}$}}
\def\ker{\mbox{$\mathop{\mbox{\rm Ker$\,$}}$}}
\def\hom{\mbox{$\mathop{\mbox{\rm Hom}}$}}

\def\o{\sp{[1]}}
\def\t{\sp{[2]}}

\nc{\spp}{\mbox{${\cal S}{\cal P}(P)$}}
\nc{\ob}{\mbox{$\Omega\sp{1}\! (\! B)$}}
\nc{\op}{\mbox{$\Omega\sp{1}\! (\! P)$}}
\nc{\oa}{\mbox{$\Omega\sp{1}\! (\! A)$}}
\nc{\dr}{\mbox{$\Delta_{R}$}}
\nc{\dsr}{\mbox{$\Delta_{\Omega^1P}$}}
\nc{\ad}{\mbox{$\mathop{\mbox{\rm Ad}}_R$}}
\nc{\as}{\mbox{$A(S^3\sb s)$}}
\nc{\bs}{\mbox{$A(S^2\sb s)$}}
\nc{\slc}{\mbox{$A(SL(2,\C))$}}
\nc{\suq}{\mbox{$\cO(SU_q(2))$}}
\nc{\tc}{\widetilde{can}}

\rnc{\epsilon}{\varepsilon}
\rnc{\phi}{\varphi}
\nc{\ha}{\mbox{$\alpha$}}
\nc{\hb}{\mbox{$\beta$}}
\nc{\hg}{\mbox{$\gamma$}}
\nc{\hd}{\mbox{$\delta$}}
\nc{\he}{\mbox{$\varepsilon$}}
\nc{\hz}{\mbox{$\zeta$}}
\nc{\hs}{\mbox{$\sigma$}}
\nc{\hk}{\mbox{$\kappa$}}
\nc{\hm}{\mbox{$\mu$}}
\nc{\hn}{\mbox{$\nu$}}
\nc{\hl}{\mbox{$\lambda$}}
\nc{\hG}{\mbox{$\Gamma$}}
\nc{\hD}{\mbox{$\Delta$}}
\nc{\hT}{\mbox{$\Theta$}}
\nc{\ho}{\mbox{$\omega$}}
\nc{\hO}{\mbox{$\Omega$}}
\nc{\hp}{\mbox{$\pi$}}
\nc{\hP}{\mbox{$\Pi$}}

\nc{\qpb}{quantum principal bundle}

\def\C{{\Bbb C}}
\def\N{{\Bbb N}}

\def\Q{{\Bbb Q}}
\def\cO{{\mathcal O}}

\def\nc{\newcommand}

\def\tc{\tilde c}

\def\ra{\rightarrow}

\def\t{\tau}
\def\te{{\tilde e}}

\def\sw#1{{\sb{(#1)}}}

\def\tens{\mathop{\otimes}}

\def\o{{}_{(1)}}
\def\t{{}_{(2)}}

\def\<{{\langle}}
\def\>{{\rangle}}

\def\id{{\rm id}}
\def\eps{\varepsilon}

\def\q2{{q^{-2}}}

\def\note#1{{}}

\def\note#1{}

\def\Rhom#1#2#3{{{\rm Hom}\sp{#1}(#2,#3)}}
\def\Lhom#1#2#3{{{}\sp{#1}{\rm Hom}(#2,#3)}}

\def\Label#1{\label{#1}\ifmmode\llap{[#1] }\else
\marginpar{\smash{\hbox{\tiny [#1]}}}\fi}
\def\Label{\label}

\def\equad{\kern -1.7em}

\def\C{{\Bbb C}}

\def\Hom{{\rm Hom}}

\begin{document}
\begin{frontmatter}
\selectlanguage{english}
\vspace*{-95pt}
\title{The Chern-Galois character}

\vspace{-2.6cm}

\selectlanguage{francais}
\title{Le caract\`ere de Chern-Galois}

\selectlanguage{english}
\author[authorlabel1]{Tomasz Brzezi\'nski},
\ead{T.Brzezinski@swansea.ac.uk}
\author[authorlabel2]{Piotr M.\ Hajac}
\ead{http://www.fuw.edu.pl/$\!\widetilde{\phantom{m}}\!$pmh}
\address[authorlabel1]{Department of Mathematics,
University of Wales Swansea,
Singleton Park, Swansea SA2 8PP, U.K.}
\address[authorlabel2]{Instytut Matematyczny, Polska Akademia Nauk,
ul.\ \'Sniadeckich 8, Warszawa, 00-956 Poland;\\
Katedra Metod Matematycznych Fizyki, Uniwersytet Warszawski
ul.\ Ho\.za 74, Warszawa, 00-682 Poland}

\begin{abstract}
Following the idea of Galois-type extensions and entwining structures,
we define the notion of a
principal extension of noncommutative algebras. We show that
modules associated to such extensions via finite-dimensional
corepresentations are finitely generated projective, and determine an
explicit formula for
the Chern character applied to the thus obtained modules.

\vskip 0.5\baselineskip

\selectlanguage{francais}
\noindent{\bf R\'esum\'e}
\vskip 0.5\baselineskip
\noindent
Nous nous inspirons des extensions de type Galois et des structures
enlac\'ees
pour d\'efinir la notion d'extension principale d'alg\`ebres
non commutatives. Nous montrons que les modules associ\'es \`a
de telles extensions au travers de corepr\'esentations de dimension
finie sont projectifs et de type fini, et nous d\'eterminons une formule
explicite pour le caract\`ere de Chern appliqu\'e aux modules ainsi
obtenus.
\end{abstract}
\end{frontmatter}

\selectlanguage{english}

\section{Introduction}
\label{}
The aim of this paper is twofold. First we need to determine a class
of Galois-type
extensions that are
sufficiently general to accommodate interesting examples and sufficiently
specific to
derive a number of desired
properties. This leads to the concept of  principal extensions. They play
the role of
algebraic analogues
of principal bundles.   To any such extension one can associate modules
much as vector bundles are associated to principal bundles.
For finite-dimensional
corepresentations these modules are always finitely generated projective
(see Theorem~\ref{thm.main}), and
thus fit the
formalism of the Chern-Connes pairing between K-theory and cyclic
cohomology \cite{c-a85}. On the other hand, just as the commutative
faithfully flat
Hopf-Galois extensions with bijective antipodes coincide with affine group
scheme
torsors of algebraic geometry, the principal extensions are
precisely  the much studied faithfully flat Hopf-Galois extensions with
bijective
antipodes whenever the defining
coaction is an algebra homomorphism (cf.\ Theorem~\ref{thm1} and \cite{ss}).
A very interesting and non-trivial example of a principal
extension encoding a noncommutative version of the instanton fibration
$SU(2)\ra S^7\ra S^4$ was recently constructed in \cite{bdct}.

The second and main outcome of this work is the construction and an explicit formula
for the Chern-Galois
character. This is a homomorphism of Abelian groups that using a principal
extension
assigns to
the isomorphism class of a finite-dimensional corepresentation the homology
class of an even cyclic cycle. This construction is in analogy with the
Chern-Weil
formalism for principal
bundles and bridges the coalgebra-Galois-extension \cite{s-hj90,bm98,bh99}
 and K-theoretic formalisms (cf.\ \cite{dgh01} for the Hopf-Galois
version). In particular, with the help of
finitely-summable Fredholm modules, it  allows one to apply
the analytic tool of the noncommutative index formula \cite{c-a85} to
compute the
$K_0$-invariants of line
bundles over generic Podle\'s spheres \cite{hms}, which are among prime
examples
 going beyond the Hopf-Galois
framework.

Except for the last formula, we work  over a general field $k$.
We use the usual
notations $\Delta c=c\o\tens c\t\in C\ot C$,
$((\id\otimes\Delta)\circ \Delta)(c) = c\sw 1\tens c\sw 2\tens c\sw 3
\in C\ot C\ot C$,
etc.,
$\Delta_V(v)=v_{(0)}\tens v\o\in V\ot C$, (summation understood)
for the  coproduct
of a coalgebra $C$, its iterations
and  a  right $C$-coaction on  $V$, respectively.
We denote the counit
of $C$ by $\eps$.
 For an algebra $A$,
 ${}_A\Hom(V,W)$ stands for the  space of left $A$-linear maps.
 Similarly,
for a coalgebra $C$, we write
$\Hom^C(V,W)$ for the space of right $C$-colinear maps.

\section{Principal extensions and strong connections}

The concept of a faithfully flat Hopf-Galois extension with a
bijective antipode
is a cornerstone of Hopf-Galois theory.
The following notion of a {\em principal extension} generalizes this
key concept
in a way that
it encompasses interesting examples escaping Hopf-Galois theory,
yet still enjoys a number of crucial properties of the aforementioned
class of
Hopf-Galois
extensions. It is an elaboration of the Galois-type extension
\cite[Definition~2.3]{bh99} (see the condition (1) below),
which evolved from
\cite[p.182]{s-hj90}, \cite{bm98} and other papers.
\\
\begin{definition}
\Label{def.principal}
Let $C$ be a coalgebra and $P$ an algebra and a right $C$-comodule via
$\Delta_P:P\ra P\ot C$. Put
$
B=P^{co C}:=\{b\in P~|~\Delta_P(bp)=b\Delta_P(p),\ \forall p\in P\}.
$
We say that the inclusion $B\inc P$ is a
{\em $C$-extension}. A $C$-extension  $B\inc P$ is called {\em principal}
 iff

(1)
$
can_R: P\ot_BP{\ra} P\ot C,\;p\ot_B p'\mapsto p\hD_P(p')
$
is bijective (Galois or freeness condition);

(2) $\psi:C\ot P{\ra} P\ot C$, $c\ot p\mapsto
can(can^{-1}(1\ot c)p)$ is bijective (invertibility of the canonical
entwining);

(3)
there is a group-like element $e\in C$ such that $\hD_P(p)=\psi(e\ot p)$,
$\forall p\in
P$ (co-augmentation);

(4)
 $P$ is $C$-equivariantly projective as a left $B$-module (existence of a
 strong
connection).
\end{definition}
~\\
The meaning of the last condition in Definition~\ref{def.principal}
is as follows.
Let $X$ be a left $B$-module and a right $C$-comodule such that the coaction is
$B$-linear.
 We say
that $X$ is a {\em $C$-equivariantly projective $B$-module}
 iff   for every $B$-linear $C$-colinear epimorphism $\pi: M\to N$ that
 is split
as a $C$-comodule map,
and for any  $B$-linear $C$-colinear homomorphism $f:X\to N$, there exists a
$B$-linear $C$-colinear  map $g: X\to M$ such that $\pi\circ g =f$.
For the trivial
$C$ we recover the
usual concept of projectivity. Much as for the trivial $C$, one can show
that
equivariant projectivity
is equivalent to the existence of a $B$-linear $C$-colinear splitting of the
multiplication map
$m:B\ot X\ra X$. If we take $A=\hom(C,k)^{op}$ to be the opposite of the
convolution
algebra of $C$,
then such a splitting is the same as a $(B,A)$-bimodule splitting of $m$.
Now one can
reverse
the argument and prove that the existence of such a $(B,A)$-bimodule
splitting is
equivalent to
$A$-equivariant projectivity defined analogously as $C$-equivariant
projectivity.
If $A$ is a commutative
ring and $B$ is an algebra over $A$, then we obtain an old concept of
relative
projectivity \cite[p.197]{ce56}).
On the other hand, as explained in \cite[p.314]{dgh01}, a
$(B,A)$-bimodule
splitting
of $m$ can be interpreted
as a Cuntz-Quillen type connection. The unitalized version of such
connections are
called strong connections.
More precisely, if $B\inc P$ is a principal $C$-extension, a
{\em strong connection}
is a unital left $B$-linear
right $C$-colinear splitting of the multiplication map
$B\ot P\ra P$
\cite[Remark~2.11]{dgh01}. The following lemma allows us to
conclude that
principal
extensions always admit
strong connections.\\
\begin{lemma}
Let $B\inc P$ be a $C$-extension satisfying conditions (1) and (3) in
Definition~\ref{def.principal}. Then $P$
is $C$-equivariantly projective as a left $B$-module {\em if and only if}
there exists a strong connection.
\end{lemma}~\\
The right-to-left part of the  assertion is immediate from the discussion
preceding
the lemma.
It is the proof of the existence of a strong connection (unital splitting)
that
requires some work. Next,
note that the conditions (2)--(3) of Definition~\ref{def.principal}
allow us to
give a symmetric formulation
of a strong connection. To begin with, one can define a left coaction
$_P\hD:P\ra C\ot P$,
$_P\hD(p)=\psi^{-1}(p\ot e)$, and prove that $can_L:P\ot_BP\ra C\ot P$,
$p\ot_Bp'\mapsto\,_P\hD(p)p'$,
is bijective. One can also show that
$can_L^{-1}\ci(\id\ot 1)=can_R^{-1}\ci(1\ot\id)$, so that the concept
of the translation map $\tau:=can_R^{-1}\ci(1\ot\id)$,
$\tau(c)=:c^{[1]}\ot_Bc^{[2]}$ (summation suppressed),
is left-right symmetric. This leads to:\\
\begin{lemma}
    Let $B\inc P$ be a principal $C$-extension, and let
$\pi_B:P\ot P\ra P\ot_B P$ be the canonical surjection.
Then the formulae
$s\mapsto (\ell:c\mapsto c^{[1]}s(c^{[2]}))$,
$\ell\mapsto(s:p\mapsto p\0\ell(p\1))$
define mutually inverse maps between the space of strong connections and
linear maps $\ell:C\ra P\ot P$ such that
$\pi_B\ci\ell=\tau$,
$(\id\ot\hD_P)\ci\ell=(\ell\ot\id)\ci\hD$,
$(_P\hD\otimes\id)\ci\ell=(\id\ot\ell)\ci\hD$,
$\ell(e)=1\ot 1$.
\Label{lemma.strong.inv}
\end{lemma}~\\
To avoid multiplying terminology,  such unital bicolinear liftings of the
translation map are also called
strong connections. Among other consequences of the principality of an
extension is
its coflatness.
Recall first that, for any right $C$-comodule $V$ with a coaction
$\Delta_V$ and a
left $C$-comodule $W$ with a coaction $_W\Delta$,  the
{\em cotensor product} is
defined as
$V\Box_{C}W := \ker (\id\ot\,_W\Delta - \hD_V\ot\id)\inc V\otimes W$.
A right (resp.\ left) $C$-comodule $M$ is said to be {\em coflat}
if the  functor  $M\Box_{C}-$
(resp.\ $-\Box_{C}M$)
is exact. Next, recall that there is a general concept of an
entwining
structure
$(A,C,\psi)$, where
$A$ is an algebra, $C$ a coalgebra, and $\psi:C\ot A\ra A\ot C$
is a
linear map
satisfying certain axioms
\cite[Definition~2.1]{bm98}. With these definitions, we obtain:\\
\begin{lemma}
Let $(A,C,{\psi})$ be an entwining structure such that
$\psi$ is bijective.
Assume also that there exists a
group-like $e\in C$ such that $A$ is a right $C$-comodule via
$\psi\ci(e\ot\id)$
and a left $C$-comodule via
$\psi^{-1}\ci(\id\ot e)$. Then
 $A$ is coflat as a right (resp.\ left) $C$-comodule if and only if
 there exists
    $j_{R}\in \Rhom C CA$ (resp.\  $j_{L}\in \Lhom C CA$) such that
    $j_{R}(e)=1$
(resp.\  $j_{L}(e)=1$).
(Here $C$ is a $C$-comodule via the coproduct.)
\Label{prop.jrl}
\end{lemma}~\\
The axioms (1)--(3) of a principal extension guarantee that
$(P,C,\psi)$ is an
entwining structure satisfying the assumptions of the
above lemma \cite[Theorem~2.7]{bh99}. Moreover,
with the help of Lemma~\ref{lemma.strong.inv}, it can be shown that
maps $j_L$ and
$j_R$ as in Lemma~\ref{prop.jrl} can be constructed for any principal
$C$-extension.
Combining together the results described in this section, one can prove the
following:\\
\begin{theorem}\label{thm1}
   Let $B\inc P$ be a principal $C$-extension. Then:

(1) There exists a strong connection.

   (2) $P$ is a projective left and right $B$-module.

   (3) $B$ is a direct summand of $P$ as a left and right $B$-module.

   (4) $P$ is a faithfully flat left and right $B$-module.

(5)   $P$ is a coflat left and right $C$-comodule.
\Label{cor.coflat}
\end{theorem}

\section{Associated projective modules and the Chern-Galois character}

 If $B\inc P$ is a principal $C$-extension
and $\phi:V_\phi\ra V_\phi\ot C$ is a finite-dimensional corepresentation,
then, using the technology from the previous section, one can produce a
short proof
that the left $B$-module ${\rm Hom}^C(V_\phi,P)$ of all
colinear maps from  $V_\phi$ to $P$ is finitely generated projective.
We call such modules associated modules, as they play the role of sections
of vector
bundles
associated to principal bundles (cf.\ \cite[Theorem~5.4]{b-t99}).
The main result of this paper is an explicit formula for an
idempotent representing an associated module.
By virtue of Theorem~\ref{cor.coflat}, we already know that there exists a
strong connection $\ell$
and a unital left $B$-linear map
$\hs_L:P\ra B$. Thus we can state:
\\
\begin{theorem}
Let $\ell$ be a strong connection on  a principal $C$-extension $B\inc P$
and
$\phi:V_\phi\ra V_\phi\ot C$ be a finite-dimensional corepresentation.
Let $\{p_\mu\}_\mu$ be a basis of $P$, $\{p^\mu\}_\mu$ its dual,
 $r_\mu:=(p^\mu\ot\id)\ci\ell$, and $\{e_{i}\}_i$ be a basis of $V_\phi$,
$\phi(e_j)=:\sum_{i=1}^{\dim V_\phi}e_{i}\ot e_{ij}$.
Take any $\hs_L\in\,_B\hom(P,B)$ such that $\hs_L(1)=1$,
and set $E_{(\mu,i)(\nu,j)}:=\hs_L(r_\mu(e_{ij})p_\nu)$,
$E:=(E_{(\mu,i)(\nu,j)})$.
Then, for some  $N\in\N$,
$E$ is a square matrix of size $N$,
$E^2=E$
and $B^NE\cong\Rhom C {V_\phi} P$ as a left $B$-module.
\label{thm.main}
\end{theorem}~\\
It is an immediate corollary of this theorem that a left $B$-module
$\Rhom C {V_\phi} P$ is always finitely generated
projective. On the other hand, take all the isomorphism classes of
finite-dimensional corepresentations of $C$
and view them as a semi-group via the direct sum. Denote by $R_f(C)$
the Grothendieck group of this
semi-group. It is now straightforward to verify that the assignment
$[\phi]\mapsto[\Rhom C {V_\phi} P]$
defines a homomorphism of Abelian groups $R_f(C)\ra K_0(B)$.
Combining this homomorphism with
the Chern character $ch_{2n}:K_0(B)\ra HC_{2n}(B)$, $n\in\N$,
(see \cite[p.264]{l-jl98}) yields a homomorphism
$chg_{2n}:R_f(C)\ra HC_{2n}(B)$. We call the collection of homomorphisms
$chg_{2n}$, $n\in\N$,
{\em the Chern-Galois character}. The main point of this work is that we
can use
Theorem~\ref{thm.main}
to determine an explicit formula for the Chern-Galois character.
To this end,
 define the character
of $\phi$ as $c_{\phi}:=\sum_{i=1}^{\dim V_\phi}e_{ii}$.
Evidently, $c_\phi$ does not depend on the choice of a basis. Next, put
$\ell(c)=:c^{<1>}\ot c^{<2>}$
 (summation understood), $c\in C$.
 Finally, for the sake of simplicity,  assume that $\Q\inc k$. Then one can
show that the Chern-Galois character comes out as
\[
\!\!\!\!\!\!\!\!\!\forall\; n\in\N:\;chg_{2n}([\phi])
=(-1)^n[{c_{\phi\, (2n+1)}}\!\!\!^{<2>}{c_{\phi\,(1)}}\!\!\!^{<1>}\ot
{c_{\phi\, (1)}}\!\!\!^{<2>}{c_{\phi\,(2)}}\!\!\!^{<1>}
\ot\cdots\ot{c_{\phi\, (2n)}}\!\!\!^{<2>}{c_{\phi\,(2n+1)}}\!\!\!^{<1>}].\;
\]
Note that, since any strong connection yields an idempotent
representing the same module and
the Chern character does not dependent on the choice of a
 representing idempotent, the
 formula for the Chern-Galois character is manifestly
independent of the choice of a strong connection appearing on the
right-hand side.\\

\noindent{\bf Acknowledgements.} \footnotesize
T.~Brzezi\'nski thanks the EPSRC for an Advanced
Research Fellowship.
P.M.~Hajac thanks the
European Commission for
the Marie Curie fellowship
 HPMF-CT-2000-00523 and the KBN for the grant 2 P03A 013 24.
 Both authors are very
grateful to  R.~Taillefer for the French translation
and to R.~Matthes for his helpful comments on the manuscript.
A preliminary version of the
full account of
this work can be found at
http://www.fuw.edu.pl/$\!\widetilde{\phantom{m}}\!$pmh.

\end{document}